\documentclass[12pt]{article}
\usepackage{mathrsfs}
\usepackage{amsfonts}
\usepackage{amssymb, amsmath,mathrsfs,amsfonts}

\numberwithin{equation}{section}

\begin{document}
\title {\large COMPACT DIFFERENCES  OF COMPOSITION OPERATORS ON HOLOMORPHIC FUNCTION SPACES IN THE UNIT
BALL\footnotetext{ {\it 2000 Mathematics Subject Classification.}
Primary: 47B38; Secondary 32A35, 32A36.} \footnotetext{  {\it Key
words and phrases.} Composition operators,  Hardy space, Bergman
spaces, compact differences.}}
\author{\normalsize  LIANGYING JIANG\footnote{Liangying Jiang is
supported by  Shanghai Education Research and Innovation Project
(No.10YZ185)
 and by Shanghai University Research Special Foundation  for Outstanding Young Teachers
 (No.sjr09015)}
 \quad CAIHENG OUYANG\footnote{ Caiheng Ouyang is
supported by the National Natural Science Foundation of China
(No.10971219)}}
\date{}
\maketitle \leftskip=10mm \rightskip=10mm \noindent{\small ABSTRACT.
We find a lower bound for the essential norm of the difference of
two composition operators acting on $H^2(B_N)$ or $A^2_s(B_N)$
($s>-1$). This result plays an important role in proving  a
necessary and sufficient condition for the difference of linear
fractional composition operators to be compact, which answers a
question posed by MacCluer and Weir in 2005.}

\leftskip=0mm \rightskip=0mm

\section*{\large  1 Introduction}

Let $B_N$ denote the open  unit ball in $\mathbb{C}^N$, with $D$ for
the unit disc  $B_1$. We write $\sigma$ to denote the normalized
Lebesgue  measure on  the unit sphere  $\partial B_N$, the Hardy
space  $H^2(B_N)$ is the set of all holomorphic functions $f$ in
$B_N$ such that $$||f||^2_2:= \sup\limits_{0<r<1}\int_{\partial
B_N}|f(r\zeta)|^2 d\sigma(\zeta)<\infty.$$ For $s>-1$, the standard
weighted Bergman space  $A^2_s(B_N)$  consists of holomorphic
functions $f$ in $B_N$ satisfying $$||f||^2_{2,s}:=
\int_{B_N}|f(z)|^2 w_s(z) d\nu(z)<\infty,$$ where
$$w_s(z)=\frac{\Gamma(N+s+1)}{\Gamma(N+1)\Gamma(s+1)}(1-|z|^2)^s$$
and $\nu$ denotes the  normalized Lebesgue volume measure on $B_N$.
Write  $H^2(B_N)=A^2_{-1}(B_N)$, it is well known that $A^2_s(B_N)$
is  a  Hilbert space of holomorphic functions in $B_N$ with the
reproducing kernel $K_z(w)=(1-<w, z>)^{-(N+s+1)}$ for any $s\ge -1$
(see \cite{Ru} and \cite{Zhu}).
\par  We consider  the composition operator $C_{\varphi}$  acting on the hardy space $H^2(B_N)$
or the Bergman spaces  $A^2_s(B_N)$ ($s> -1$), defined by
$C_{\varphi}f=f\circ\varphi,$ where $\varphi$  is an  analytic map
from $B_N$ into $B_N$. When $N=1$, the Littlewood Subordination
Theorem  shows that $C_{\varphi}$ is bounded for any analytic
self-map  $\varphi$  of $D$, and many other properties of
$C_{\varphi}$ have been characterized, see the good monographs
\cite{S} and  \cite{CM1} for details. However, for $N\ge 2$, one may
find examples of $\varphi:B_N\rightarrow B_N$ such that
$C_{\varphi}$ is not bounded (see Section 3.5 in \cite{CM1}).
Moreover, some  basic properties of composition operators in the
setting of the  ball are not easily managed.  The purpose of this
paper  is to characterize those pairs $\varphi$ and $\psi$ for which
the difference $C_{\varphi}-C_{\psi}$ is compact acting on
$H^2(B_N)$ or $A^2_s(B_N)$  ($s>-1$). From these results, one may
derive some information about the structure of the space of
composition operators.
\par  The topological structure of the set of composition operators
on $H^2(D)$  was first studied by Berkson in \cite{Be}. Shapiro and
Sundberg \cite{SS} improved the result of Berkson and raised the
problem   about compact differences of composition operators. In
\cite{SS} they found a lower bound for the essential norm
$||C_{\varphi}-C_{\psi}||_e$ in terms of the measure of the set
$E_{\varphi}=\{\zeta\in \partial D: |\varphi(\zeta)|=1\}$, where
$\varphi$ and $\psi$ are analytic self-maps of  $D$  and
$\varphi(\zeta):=\lim\limits_{r\to 1}\varphi(r\zeta)$. This result
has been extended to the case of  Hardy spaces $H^p(B_N)$ $(0<p\le
\infty)$   (see \cite{He} and \cite{HM}). On the other hand, using
the angular derivative, MacCluer \cite{M1} discussed the differences
of composition operators  and obtained the following essential norm
estimate
$$||C_{\varphi}-C_{\psi}||^2_e\ge |\varphi'(\zeta)|^{-\beta}$$
with $\beta=1$ for    $H^2(D)$  and $\beta=s+2$ for $A^2_s(D)$
($s>-1$), where $\varphi'(\zeta)$ is the angular derivative of
$\varphi$ at $\zeta\in \partial D$.  Thus, from the results of
Shapiro and Sundberg \cite{SS} and MacCluer \cite{M1}, one may
determine for which pairs $\varphi$ and $\psi$ the difference
$C_{\varphi}-C_{\psi}$ is compact. Recently, Aleksandrov-Clark
measures   also  have been used to study the compactness of
differences and linear combinations of composition operators on the
spaces mentioned (see  \cite{GGMS}, \cite{KM}, \cite{S98}).

In   Section 2 of this paper, we would expect  similar results about
the compact differences of composition operators  on  $H^2(D)$ and
$A^2_s(D)$ ($s>-1$) to hold in several variables.  First,  motivated
by the work of MacCluer  \cite{M1}, we will find  a lower bound for
the essential norm of composition operator difference
$C_{\varphi}-C_{\psi}$ on $H^2(B_N)$ or  $A^2_s(B_N)$ ($s>-1$), but
with
$$d_{\varphi}(\zeta)=\liminf\limits_{z\to
\zeta}\frac{1-|\varphi(z)|}{1-|z|}$$ instead of $|\varphi'(\zeta)|$
(see Theorem 2.1 in Section 2). In fact,  the Julia-Carath\'{e}odory
theorem in the disc shows that if $\varphi$ has finite angular
derivative at $\zeta\in \partial D$ then
$|\varphi'(\zeta)|=d_{\varphi}(\zeta)$. So the
Julia-Carath\'{e}odory theory in $B_N$ (see \cite{Ru} or \cite{CM1})
will be  a key  tool for its proof. Note that in the proof of
MacCluer's result (Theorem 2.2 of \cite{M1}), the main idea is to
use
$$\lim\limits_{z\to
\zeta}\frac{1-|\varphi(z)|}{1-|z|}=|\varphi'(\zeta)|$$
 as $z$ approaches $\zeta\in\partial D$ nontangentially, which is a result of the
 Julia-Carath\'{e}odory Theorem in the disc (see \cite{CM1}).
 However, for  higher dimensions, we have not found  the  corresponding result
 and some   techniques will be needed.  Moreover,  our
method can be generalized to estimate essential norms of linear
combinations of composition operators. As a consequence, we obtain
some necessary conditions for differences or linear combinations of
composition operators to be compact on $H^2(B_N)$ or $A^2_s(B_N)$
($s>-1$).
\par In  another direction, Bourdon \cite{Bo} treated the question
on compact differences    of  linear fractional composition
operators and proved that $C_{\varphi}-C_{\psi}$  is  compact  on
$H^2(D)$ if and only if both $C_{\varphi}$ and $C_{\psi}$  are
compact or $\varphi=\psi$. This result  also holds on $A^2_s(D)$
($s>-1$) from Moorhouse's result \cite{Mo}. For the linear
fractional self-map $\varphi$  of $B_N$ with a boundary fixed point,
MacCluer and Weir \cite{MW} showed that the difference
$C_{\varphi\circ\sigma}-C_{\sigma\circ\varphi}$ is compact on
$H^2(B_N)$ or   $A^2_s(B_N)$  ($s>-1$) if and only if
$\varphi\circ\sigma=\sigma\circ\varphi$, where $\sigma$  is the
adjoint map of $\varphi$,   and  asked the following question:
\\ \\ ($\ast$) \ {\it  For distinct linear fractional self-maps $\varphi$
and $\psi$ of $B_N$ can $C_{\varphi}-C_{\psi}$ ever be compact?}
\\ \par In Section 3, we  then focus on compact differences
of linear fractional composition operators. For linear fractional
self-maps $\varphi$ and $\psi$ of $B_N$, we will prove that
$C_{\varphi}-C_{\psi}$ is compact on $H^2(B_N)$ or  $A^2_s(B_N)$
($s>-1$) if and only if both $C_{\varphi}$ and $C_{\psi}$  are
compact or $\varphi=\psi$ (Theorem 3.1 in Section 3), which  answers
the question ($\ast$). The basic ideas  come from   Bourdon
\cite{Bo} and MacCluer and Weir \cite{MW}. In our proof an important
tool is the result about   compact difference of composition
operators in Section 2, which says that if $C_{\varphi}-C_{\psi}$ is
compact then $\varphi(\zeta)=\psi(\zeta)$ and
$d_{\varphi}(\zeta)=d_{\psi}(\zeta)$ at some point $\zeta\in\partial
B_N$ ( see Corollary 2.2 in Section 2).  This will give  a very
useful information for  the relations of  matrixes associated with
$\varphi$ and $\psi$. In this point, our approach  is different from
that used by MacCluer and Weir \cite{MW}.

For the proof of MacCluer and Weir's result \cite{MW}, under the
condition of $\varphi$ having a boundary fixed point, assume $e_1$
to be fixed, they obtained
$\varphi\circ\sigma(e_1)=\sigma\circ\varphi(e_1)$.   They also found
that the adjoint maps of $\varphi\circ\sigma$ and
$\sigma\circ\varphi$ are  themselves and then  deduced that
 $D_1(\varphi\circ\sigma)_1(e_1)=D_1(\sigma\circ\varphi)_1(e_1)=1$.
According to the proof of Lemma 4.2 in \cite{JO}, we see that this
result always holds  in the case of $\varphi$ fixing
 $e_1$. Note that the
Julia-Carath\'{e}odory Theorem in $B_N$ gives
$D_1(\varphi\circ\sigma)_1(e_1)=d_{\varphi\circ\sigma}(e_1)$ and
$D_1(\sigma\circ\varphi)_1(e_1)=d_{\sigma\circ\varphi}(e_1)$. Hence,
if  $\varphi$ fixes  a boundary point,  MacCluer and Weir in fact
obtained the same result as ours, that is
$\varphi\circ\sigma(e_1)=\sigma\circ\varphi(e_1)$ and
$d_{\varphi\circ\sigma}(e_1)=d_{\sigma\circ\varphi}(e_1)$.  However,
if  $||\varphi||_\infty=1$,  this will hold automatically from the
compactness of $C_{\varphi\circ\sigma}-C_{\sigma\circ\varphi}$  by
Corollary 2.2 in Section 2. Thus, the hypothesis  that   $\varphi$
fixes  a boundary point  in their result can be replaced by a weaker
condition  $||\varphi||_\infty=1$ (see Theorem 3.2 in Section 3).

This work is part of the first author's doctoral thesis (see
\cite{JO1}), but, at that time,  the method for proving
$\gamma_k=\gamma'_k$ ($k=1,\cdots,n$) in the proof of  Theorem 3.1
was not correct. In this paper, we improve the proof of  Theorem 3.1
 and obtain  some other results. Recently, the authors
learned that Heller et al \cite{HMW} independently proved Theorem
3.1  using different methods.

\section*{\large 2 Essential norms of composition operator
differences and linear combinations}

The essential norm of an operator $T$ on the space $\mathcal {H}$ is
defined by $||T||_e=\inf\{||T-K||: K  \  \mbox{is compact on }
\mathcal {H}\}$. In \cite{M1}, MacCluer considered the topological
space of composition operators and obtained  the following result.
\\ \\ {\bf Theorem A.} \begin{em}  Let $\varphi,\psi: D\rightarrow D$ be
analytic maps and suppose that $\varphi$ has a finite angular
derivative at $\zeta\in \partial D$.  Consider $C_{\varphi}$ and
$C_{\psi}$ acting on $H^2(D)$ or  $A^2_s(D)$ for $s>-1$. Then unless
$\psi(\zeta)=\varphi(\zeta)$ and $\psi'(\zeta)=\varphi'(\zeta)$, one
have $$||C_{\varphi}-C_{\psi}||^2_e\ge |\varphi'(\zeta)|^{-\beta},$$
where $\beta=1$ for the space $H^2(D)$  and $\beta=s+2$ for  the
spaces $A^2_s(D)$. \end{em}
\\ \par If $\varphi$ and $\psi$ have radial limits of modulus $1$ at  $\zeta\in \partial
D$ with  $\varphi(\zeta)=\psi(\zeta)$ and
$|\varphi'(\zeta)|=|\psi'(\zeta)|$, we say that  $\varphi$ and
$\psi$ have the same data at this point (see \cite{M1}).
Immediately, from Theorem A, one may  get that if
$C_{\varphi}-C_{\psi}$ is compact then $\varphi$ and $\psi$ must
have the same data for  those points, at which $\varphi$ has finite
angular derivatives.
\par  In this section, we will discuss the analogue of  Theorem A
for the ball, but  in higher dimensions, our lower bound needs  a
corresponding form of
 the  angular derivative $|\varphi'(\zeta)|$  of  $\varphi:
D\rightarrow D$.  First, we summarize some relevant results on the
angular derivative and the Julia-Carath\'{e}odory theory in the
ball.
\par A curve $\Gamma$ in $B_N$ will be called a $\zeta$-- curve if $\Gamma$
approaches a point $\zeta\in \partial B_N$. We say that a function
$f:  B_N\rightarrow \mathbb{C}$ has restricted limit $L$ at
$\zeta\in \partial B_N$, if $\lim\limits_{t\to 1}f(\Gamma(t))=L$ for
every $\zeta$-- curve $\Gamma(t)$ that satisfies
$$\lim\limits_{t\to 1} \frac{|\Gamma(t)-\gamma(t)|^2}{1-|\gamma(t)|^2}=0$$
and $$\frac{|\gamma(t)-\zeta|}{1-|\gamma(t)|}\le M<\infty \ \ \
\quad  \mbox{for} \  0\le t<1,$$ where $\gamma(t)=<\Gamma(t),
\zeta>\zeta$ is the  projection of $\Gamma$ onto the complex line
through $\zeta$. In this case, the curve $\Gamma$ is said to be
restricted and its orthogonal projection $\gamma$ is nontangential
(see \cite{Ru}).
\par Let $\varphi$ be an analytic self-map of $B_N$  and $\zeta\in \partial B_N$,
if there exists  a point  $\eta\in \partial B_N$ such that the
restricted limit of $$\frac{<\eta-\varphi(z), \eta>}{ <\zeta-z,
\zeta>}$$  exists then $\varphi$ is said to have finite angular
derivative at $\zeta$. By the Julia-Carath\'{e}odory  Theorem in
$B_N$, this is equivalent to
$$d_{\varphi}(\zeta)=\liminf\limits_{z\to
\zeta}\frac{1-|\varphi(z)|}{1-|z|}<\infty,$$  where $z$ approaches
$\zeta$  unrestrictedly in $B_N$. Moreover, under these conditions,
$\varphi$ has restricted limit $\eta$ at $\zeta$  and
$D_{\zeta}\varphi_{\eta}(z)=<\varphi'(z)\zeta, \eta>$  has
restricted limit $d_{\varphi}(\zeta)$.
\par Next, making use of the Julia-Carath\'{e}odory  Theorem in
$B_N$,  we will give lower bounds for essential norms of differences
and linear combinations of composition operators on $H^2(B_N)$ or
$A^2_s(B_N)$ ($s>-1$). Therefore,   some information about the
compactness of them can be obtained. \\ \\ {\bf Theorem
2.1.}\begin{em} Let $\varphi$ and $\psi$ be  analytic self-maps of
$B_N$. Suppose that they induce bounded composition operators    on
$H^2(B_N)$ or $A^2_s(B_N)$ ($s>-1$) and $\varphi$ has finite angular
derivative at $\zeta\in\partial B_N$. Then, unless
$\psi(\zeta)=\varphi(\zeta)$ (as radial limits) and
$d_{\psi}(\zeta)=d_{\varphi}(\zeta)$, we have
$$||C_{\varphi}-C_{\psi}||^2_e\ge d_{\varphi}(\zeta)^{-\beta},$$
where $\beta=N$ for $H^2(B_N)$  and $\beta=N+s+1$ for
$A^2_s(B_N)$.\end{em} \\ \\  Proof. If $\varphi$ has finite angular
derivative at $\zeta\in\partial B_N$, by the Julia-Carath\'{e}odory
Theorem in $B_N$, there exists a point $\eta\in \partial B_N$ such
that $\varphi(\zeta):=\lim_{r\to 1}\varphi(r\zeta)=\eta$. Assume
that $U$ and $V$ are unitary transformations on $B_N$ which send
$e_1$ to $\zeta$ and $\eta$ respectively, where
$e_1=(1,0,\ldots,0)=(1, 0')$. Let $V^{\ast}$ be the adjoint of $V$
with  $V^{\ast}=V^{-1}$. Then the map $\phi(z)=V^{\ast}\varphi U(z)$
also has finite angular derivative at $e_1$ and
{\setlength\arraycolsep{2pt}
\begin{eqnarray*}d_{\phi}(e_1)&=&\lim_{r\to 1}D_1 \phi_1(r e_1)=\lim_{r\to 1}<\phi'(r
e_1)e_1, e_1> =\lim_{r\to 1}<V^{\ast}\varphi'(U(r e_1)) U e_1,
e_1>\\ &=&\lim_{r\to 1}<\varphi'(r \zeta)\zeta, \eta>=\lim_{r\to
1}D_{\zeta}\varphi_{\eta}(r\zeta)=
d_{\varphi}(\zeta),\end{eqnarray*}}where $\phi_1$ denotes the first
coordinate function of $\phi$.  Moreover, write $\tau=V^{\ast}\psi
U$, we have
$$||C_{\phi}-C_{\tau}||_e=||C_{V^{\ast}\varphi U}-C_{V^{\ast}\psi
U}||_e=||C_{U}(C_{\varphi}-C_{\psi})C_{V^*}||_e=||C_{\varphi}-C_{\psi}||_e.
$$
So the proof  will be complete if we  can show that the result holds
for  $\phi$ and $\tau$. Thus, we may assume $\zeta=\eta=e_1$.
\par Let $K_z$ be the reproducing kernel  for $z\in B_N$,  it is
easy to see
$$||C_{\varphi}-C_{\psi}||_e^2=||(C_{\varphi}-C_{\psi})^{\ast}||_e^2\ge
\limsup\limits_{|z|\to 1}\frac{||(C_{\varphi}-C_{\psi})^{\ast}
K_z||^2}{||K_z||^2}.$$ Since $C_{\varphi}^{\ast}K_z=K_{\varphi(z)}$,
we can  write {\setlength\arraycolsep{2pt}
\begin{eqnarray*}\frac{||(C_{\varphi}-C_{\psi})^{\ast}
K_z||^2}{||K_z||^2}
&=&\frac{||K_{\varphi(z)}-K_{\psi(z)}||^2}{||K_z||^2}\\
&=&(1-|z|^2)^{\beta} [||K_{\varphi(z)}||^2+||K_{\psi(z)}||^2-
2\mbox{Re}\, K_{\varphi(z)}(\psi(z))] \\ &=&
\biggl(\frac{1-|z|^2}{1-|\varphi(z)|^2}\biggr)^{\beta}+\biggl(\frac{1-|z|^2}
{1-|\psi(z)|^2}\biggr)^{\beta}
-2\mbox{Re}\biggl(\frac{1-|z|^2}{1-<\psi(z),
\varphi(z)>}\biggr)^{\beta},
\end{eqnarray*}}where the norm $||\cdot||$  is in the space $H^2(B_N)$
or the spaces  $A^2_s(B_N)$ ($s>-1$), and $K_z(w)=(1-<w,
z>)^{-\beta}$ is the corresponding  reproducing kernel (see Section
1).
\par Our goal is to estimate the first term and the third term on the last line of the equation above. If
$\psi(e_1):=\lim\limits_{r\to 1} \psi(r e_1)\ne e_1$, then there
exists a sequence $\{r_n\}$  going  to 1 as $n\to\infty$ such that
$\lim\limits_{n\to \infty} \psi (r_n e_1)=w\ne e_1$,  which implies
$$\lim_{n\to \infty}\mbox{Re}\biggl(\frac{1-r_n^2}{1-<\psi(r_n e_1),
\varphi(r_n e_1)>}\biggr)^{\beta}=0.$$  On the  other hand, the
proof of the Julia-Carath\'{e}odory  Theorem in $B_N$ gives
$$\lim_{n\to \infty}\frac{1-|\varphi(r_n
e_1)|}{1-r_n}=d_{\varphi}(e_1).$$ Thus,  we deduce that
$||C_{\varphi}-C_{\psi}||_e^2\ge d_{\varphi}(e_1)^{-\beta}$.
\par Next, if  $\lim\limits_{r\to 1} \psi(r e_1)=e_1$ but $d_{\psi}(e_1)\ne
d_{\varphi}(e_1)$, to deal with this case, the argument  used to
prove Theorem A is not helpful.  For a self-map $\varphi$ of $D$, if
$\varphi$ has finite angular derivative at $\zeta\in\partial D$,
then
$$\lim\limits_{z\to \zeta}\frac{1-|\varphi(z)|}{1-|z|}$$ has
nontangential limit $|\varphi'(\zeta)|$  by  the
Julia-Carath\'{e}odory  theory in the disc. However, we  don't know
whether this  would happen  in the setting of the ball and we need
some different approaches.  For any  curve $\gamma$ approaching $1$
nontangentially in $D$, the curve $\Gamma\equiv\{z=(\lambda, 0'):
\lambda\in \gamma\}$ is a restricted $e_1$-- curve in $B_N$. Note
that  the angular derivative of $\varphi$  existing implies that
$(1-\varphi_1(z))/(1-z_1)$ has restricted limit $d_{\varphi}(e_1)$
at $e_1$,  so it tends to $d_{\varphi}(e_1)$  as $z$ approaches
$e_1$ along $\Gamma$.  Define $\rho(\lambda)=\varphi_1(\lambda
e_1)=\varphi_1(\lambda, 0')$ for $\lambda\in D$,  the  above
discussion shows that
$$\frac{1-\rho(\lambda)}{1-\lambda}=\frac{1-\varphi_1(\lambda, 0')}{1-\lambda}$$
has finite nontangential limit $d_{\varphi}(e_1)$ as $\lambda\to 1$.
Therefore,   the map $\rho$ has finite angular derivative at $1$. By
the Julia-Carath\'{e}odory  Theorem in $D$,  we see that
$$\lim_{\lambda\to
1}\frac{1-|\rho(\lambda)|}{1-|\lambda|}=d_{\varphi}(e_1)$$ as
$\lambda$ approaches $1$ nontangentially. Combining
$$\frac{1-|\varphi(z_1, 0')|}{1-|z_1|}\le \frac{1-|\varphi_1(z_1, 0')|}{1-|z_1|}
=\frac{1-|\rho(z_1)|}{1-|z_1|}$$ with $\liminf\limits_{z\to
e_1}\frac{1-|\varphi(z)|}{1-|z|}=d_{\varphi}(e_1)$, we get
$$\lim_{z_1\to 1}\frac{1-|\varphi(z_1,
0')|}{1-|z_1|}=d_{\varphi}(e_1)$$ as $z_1$ tends to $1$
nontangentially.
\par Now, we discuss two cases for $d_{\psi}(e_1)\ne
d_{\varphi}(e_1)$. First, if $d_{\psi}(e_1)<\infty$, we compute that
{\setlength\arraycolsep{2pt}
\begin{eqnarray*} \frac{1-<\psi(z), \varphi(z)>}{1-|z|^2}
&=&\frac{1-|\varphi(z)|^2}{1-|z|^2}+\frac{<\varphi(z)-\psi(z),
\varphi(z)>}{1-|z|^2}\\ &=&
\frac{1-|\varphi(z)|^2}{1-|z|^2}+\frac{1-z_1}{1-|z|^2}\biggl[\biggl(\frac{1-\psi_1(z)}{1-z_1}
-\frac{1-\varphi_1(z)}{1-z_1}\biggr)\overline{\varphi_1(z)}\\
&&+\sum_{j=2}^N\biggl(\frac{|\varphi_j(z)|^2}
{1-z_1}-\frac{\psi_j(z)\overline{\varphi_j(z)}}{1-z_1}\biggr)\biggr].\end{eqnarray*}}Let
$\Gamma_{e_1,M}\equiv\{z=(z_1, 0')\in B_N:
\frac{|1-z_1|}{1-|z_1|^2}=M\}$. It is clear that  $\Gamma_{e_1,M}$
is  a restricted $e_1$-- curve  and the orthogonal projection of
$\Gamma_{e_1,M}$ is nontangengtial.  As $z$  approaches $e_1$ along
$\Gamma_{e_1,M}$, the previous argument shows that
$$\lim\limits_{z\to
e_1}\frac{1-|\varphi(z)|^2}{1-|z|^2}=d_{\varphi}(e_1).$$ On the
other hand, by the Julia-Carath\'{e}odory  theory  in $B_N$, we have
$$  \frac{1-\varphi_1(z)}{1-z_1}\to d_{\varphi}(e_1),
 \qquad  \frac{1-\psi_1(z)}{1-z_1}\to d_{\psi}(e_1),$$
$$ \frac{\varphi_j(z)}{(1-z_1)^{1/2}}\to 0,
 \qquad \frac{\psi_j(z)}{(1-z_1)^{1/2}}\to 0  \qquad  \mbox{for} \ \ 2\le j \le N$$
and $\varphi_1(z)$ has finite limit $1$  as $z\to e_1$  along
$\Gamma_{e_1,M}$. Note that $|(1-\overline{z_1})/(1-z_1)|=1$, so
$$\frac{\overline{\varphi_j(z)}}{(1-z_1)^{1/2}}=
\overline{\biggl(\frac{\varphi_j(z)}{(1-z_1)^{1/2}}\biggr)}\cdot\biggl(\frac{1-\overline{z_1}}{1-z_1}\biggr)^{1/2}
\to 0$$ holds for $2\le j \le N$.  Write $1-z_1=|1-z_1|e^{i\theta}$,
all these results yield that
$$\lim\limits_{\substack{z\in \Gamma_{e_1,M} \\  z\to
e_1}}\frac{1-<\psi(z),
\varphi(z)>}{1-|z|^2}=d_{\varphi}(e_1)+Me^{i\theta}(d_{\psi}(e_1)-d_{\varphi}(e_1)).$$
It follows that
$$\lim\limits_{\substack{z\in \Gamma_{e_1,M} \\  z\to
e_1}}\mbox{Re}\biggl(\frac{1-|z|^2}{1-<\psi(z),
\varphi(z)>}\biggr)^{\beta}$$  converges to $0$ as $M$ tends to
infinity. Consequently, we obtain $||C_{\varphi}-C_{\psi}||^2_e\ge
d_{\varphi}(e_1)^{-\beta}$.

If $d_{\psi}(e_1)=\infty$, we use the following inequality (see
\cite{MW}) $$\frac{||(C_{\varphi}-C_{\psi})^{\ast}
K_z||^2}{||K_z||^2}\ge
\biggr(\frac{||K_{\varphi(z)}||-||K_{\psi(z)}||}{||K_z||}\biggl)^2+2(1-u(z))
\frac{||K_{\varphi(z)}||}{||K_z||}\cdot\frac{||K_{\psi(z)}||}{||K_z||},$$
where $u(z)=(1-\rho^2(z))^{\beta/2}$ and $\rho(z)$ is the
pseudohyperbolic distance between $\varphi(z)$ and $\psi(z)$  with
$$1-\rho^2(z)=\frac{(1-|\varphi(z)|^2)(1-|\psi(z)|^2)}{|1-<\varphi(z),
\psi(z)>|^2}.$$  Note that the second term on the right side of the
inequality is not less than $0$  and
$$\lim\limits_{r\to 1}\frac{||K_{\varphi(r e_1)}||}{||K_{r e_1}||}=\lim\limits_{r\to 1}
\biggl(\frac{1-r^2} {1-|\varphi(re_1)|^2}\biggr)^{\frac{\beta}{2}}=
d_{\varphi}(e_1)^{-\frac{\beta}{2}}.$$ At the same time,
$d_{\psi}(e_1)=\liminf\limits_{z\to
 e_1}\frac{1-|\psi(z)|^2}{1-|z|^2}=\infty$ implies
 $$\lim_{r\to 1}\frac{||K_{\psi(re_1)}||}{||K_{re_1}||}=
 \lim_{r\to 1}\biggl(\frac{1-r^2}{1-|\psi(re_1)|^2}\biggr)^{\frac{\beta}{2}}=0.$$
Thus, we have
$$||C_{\varphi}-C_{\psi}||^2_e\ge\lim_{r\to 1}\frac{||(C_{\varphi}-C_{\psi})^{\ast}
K_{re_1}||^2}{||K_{re_1}||^2}\ge\lim\limits_{r\to 1}
\biggr(\frac{||K_{\varphi(re_1)}||-||K_{\psi(re_1)}||}{||K_{re_1}||}\biggl)^2
=d_{\varphi}(e_1)^{-\beta}$$  as desired.  \ \ $\Box$

As a corollary of Theorem 2.1, we get a necessary condition for  the
difference $C_{\varphi}-C_{\psi}$ to be compact. This result will
provide some heuristics for the proof of our theorem  in Section 3.
\\ \\ {\bf Corollary 2.2.}\begin{em} Suppose that $C_{\varphi}$ and
$C_{\psi}$ are bounded  on  $H^2(B_N)$ or    $A^2_s(B_N)$ ($s>-1$).
If
  $C_{\varphi}-C_{\psi}$ is compact  then $\varphi(\zeta)=\psi(\zeta)$
and $d_{\varphi}(\zeta)=d_{\psi}(\zeta)$ must hold at the  point
$\zeta\in \partial B_N$, where the angular derivative of  $\varphi$
exists.\end{em}
\\ \par In fact, for the case $d_{\psi}(e_1)=\infty$
in the proof of Theorem 2.1,  using similar idea of Kriete and
Moorhouse in \cite{KM}, we have $$
|K_{\varphi(z)}(\psi(z))|=|<K_{\varphi(z)}, K_{\psi(z)}>| \le
||K_{\varphi(z)}||^{1/2}||K_{\psi(z)}||^{1/2}$$ from the Schwarz
inequality. That is
$$\biggl(\frac{1-|z|^2}{|1-<\psi(z),
\varphi(z)>|}\biggr)^{\beta}\le
\biggl(\frac{1-|z|^2}{1-|\varphi(z)|^2}\biggr)^{\beta/2}
\biggl(\frac{1-|z|^2}{1-|\psi(z)|^2}\biggr)^{\beta/2}.$$ In the
proof of Theorem 2.1, we obtain
$$\lim\limits_{\substack{z\in \Gamma_{e_1,M} \\  z\to
e_1}}\frac{1-|z|^2}{1-|\varphi(z)|^2}=\frac{1}{d_{\varphi}(e_1)}.$$
Now, $(1-|z|^2)/(1-|\psi(z)|^2)$ tends to $0$  as $z\to e_1$
unrestrictedly.  Hence,  in the case $d_{\psi}(e_1)=\infty$,
$$\lim\limits_{\substack{z\in \Gamma_{e_1,M} \\  z\to
e_1}}\frac{1-|z|^2}{|1-<\psi(z), \varphi(z)>|}=0$$ holds and  we
also get that  $||C_{\varphi}-C_{\psi}||^2_e\ge
d_{\varphi}(e_1)^{-\beta}$. Now, combining the above  discussion
with the proof of Theorem 2.1, for analytic self-maps $\varphi$ and
$\psi$ of $B_N$ with $\varphi(e_1)=e_1$, we deduce that
\begin{eqnarray*}\lim\limits_{M\to\infty}\lim\limits_{\substack{z\in \Gamma_{e_1,M} \\  z\to
e_1}}\frac{1-|z|^2}{1-<\psi(z), \varphi(z)>}=
\begin{cases}\frac{1}{d_\varphi(e_1)}, &
\begin{subarray}{l}\mbox{if} \ \psi(e_1)=\varphi(e_1) \\ \mbox{and}\
d_{\psi}(e_1)=d_{\varphi}(e_1)<\infty,\end{subarray}
\\ 0, & \mbox{otherwise.}
\end{cases}
\end{eqnarray*} Therefore, we have the following result
for linear combinations of composition operators.
\\ \\ {\bf Theorem 2.3.}\begin{em} Let $\phi_1,\ldots, \phi_m$  be a
class of analytic self-maps of $B_N$, which induce bounded
composition operators on $H^2(B_N)$ or    $A^2_s(B_N)$ ($s>-1$).
Then for any complex numbers $c_1,\ldots, c_m$ and $\zeta\in
\partial B_N$,
$$||c_1\phi_1+\cdots+c_m\phi_m||^2_e\ge
\sum\limits_{d_{\phi_j}(\zeta)<\infty}\biggl|\sum\limits_{\substack{\phi_l(\zeta)=\phi_j(\zeta)\\
d_{\phi_l}(\zeta)=d_{\phi_j}(\zeta)}}
c_l\biggr|^2\frac{1}{d_{\phi_j}(\zeta)^\beta}
$$
with $\beta=N$ for $H^2(B_N)$  and $\beta=N+s+1$ for
$A^2_s(B_N)$.\end{em} \\ \\   Proof.  As in the proof of Theorem
2.1, we may assume $\zeta=e_1$. It is clear that
{\setlength\arraycolsep{2pt}\begin{eqnarray*}
||c_1\phi_1+\cdots+c_m\phi_m||^2_e&\ge&\limsup\limits_{|z|\to
1}\biggl|\biggl|\biggl(c_1\phi_1+\cdots+c_m\phi_m\biggr)^*\frac{K_z}{||K_z||}\biggr|\biggr|^2
\\ &\ge&
\lim\limits_{M\to\infty}\lim\limits_{\substack{z\in \Gamma_{e_1,M}
\\  z\to
e_1}}\biggl|\biggl|\biggl(\overline{c_1}\phi_1^*+\cdots+\overline{c_m}\phi_m^*\biggr)\frac{K_z}
{||K_z||}\biggr|\biggr|^2
\\ &=&\sum\limits_{j,l=1}^m \overline{c_j}\, c_l \lim\limits_{M\to\infty}\lim\limits_{\substack{z\in \Gamma_{e_1,M}
\\  z\to
e_1}}\biggl(\frac{1-|z|^2}{1-<\phi_l(z), \phi_j(z)>}\biggr)^{\beta}.
\end{eqnarray*}}If  $\phi_j$  has finite angular derivative at $e_1$ and
$\phi_j(e_1)\ne e_1$ for some $j$, there exists a unitary
transformation  $W$ such that $W\phi_j(e_1)=e_1$, then
$$\frac{1-|z|^2}{1-<W\phi_l(z), W\phi_j(z)>}=\frac{1-|z|^2}{1-<\phi_l(z), \phi_j(z)>}$$
and $d_{W\phi_j}(e_1)=d_{\phi_j}(e_1)$. Since the result proceeding
Theorem 2.3 holds for the map $W\phi_j$, we  then obtain that
\begin{eqnarray*}\lim\limits_{M\to\infty}\lim\limits_{\substack{z\in \Gamma_{e_1,M} \\  z\to
e_1}}\frac{1-|z|^2}{1-<\phi_l(z), \phi_j(z)>}=
\begin{cases}\frac{1}{d_{\phi_j}(e_1)}, &
\begin{subarray}{l}\mbox{if} \ \phi_l(e_1)=\phi_j(e_1) \\ \mbox{and}\
d_{\phi_l}(e_1)=d_{\phi_j}(e_1)<\infty,\end{subarray}
\\ 0, & \mbox{otherwise.}
\end{cases}
\end{eqnarray*}Therefore,
{\setlength\arraycolsep{2pt}\begin{eqnarray*}
||c_1\phi_1+\cdots+c_m\phi_m||^2_e&\ge&
\sum\limits_{j=1}^m\sum\limits_{l=1}^m \overline{c_j}\, c_l
\lim\limits_{M\to\infty}\lim\limits_{\substack{z\in \Gamma_{e_1,M}
\\  z\to e_1}}\biggl(\frac{1-|z|^2}{1-<\phi_l(z),
\phi_j(z)>}\biggr)^{\beta}
\\
&=&\sum\limits_{j=1}^m \biggl(\sum\limits_{\substack{\phi_l(e_1)=\phi_j(e_1)\\
d_{\phi_l}(e_1)=d_{\phi_j}(e_1)<\infty}}\overline{c_j}\,
c_l\frac{1}{d_{\phi_j}(e_1)^\beta}\biggr)
\\
&=&\sum\limits_{d_{\phi_j}(e_1)<\infty}\biggl|\sum\limits_{\substack{\phi_l(e_1)=\phi_j(e_1)\\
d_{\phi_l}(e_1)=d_{\phi_j}(e_1)}}
c_l\biggr|^2\frac{1}{d_{\phi_j}(e_1)^\beta},\end{eqnarray*}} which
is the desired conclusion. \ \ $\Box$
\par Immediately, we have the following result  for the compactness
of linear fractional combinations of composition operators.
\\ \\ {\bf Corollary 2.4.}\begin{em} Suppose that
$C_{\phi_1},\ldots,C_{\phi_m}$ are bounded and
$c_1\phi_1+\cdots+c_m\phi_m$ is compact when acting on $H^2(B_N)$ or
$A^2_s(B_N)$ ($s>-1$). For any point $\zeta\in\partial B_N$ at which
if $\phi_j$ has finite angular derivative for some $j=1,\ldots,m$,
then
$$\sum\limits_{\substack{\phi_l(\zeta)=\phi_j(\zeta)\\
d_{\phi_l}(\zeta)=d_{\phi_j}(\zeta)}} c_l=0.$$\end{em}

\section*{\large 3 Compact differences of linear fractional
composition operators}

A linear fractional map of $\mathbb{C}^N$ is defined by
$$\varphi(z)=\frac{Az+B}{<z, C>+d},$$ where $A=(a_{jk})$ is  an $N\times N$ matrix,
$B=(b_j)$, $C=(c_j)$ are $N\times 1$ column vectors, and $d$ is   a
complex number.  The matrix
\begin{displaymath} m_{\varphi}=\left(\begin{array}{ccc}
 A & B \\ C^{\ast} & d \end{array}\right) \end{displaymath}is
 called  a matrix associated with $\varphi$ and we  write
$\varphi\sim m_{\varphi}$. If  $\varphi$ is a linear fractional map
from $B_N$ into $B_N$, Cowen and MacCluer \cite{CM2} proved that the
adjoint map $\sigma$  given by
$$\sigma(z)=\frac{A^*z-C}{<z, -B>+\overline{d}}$$  maps $B_N$
into itself. In \cite{CM2} they   also  deduced that the composition
operator $C_\varphi$ is bounded on $H^2(B_N)$   and $A^2_s(B_N)$
($s>-1$).

In \cite{MW}, MacCluer and Weir considered the essential normality
of linear fractional composition operators   on $H^2(B_N)$ or
$A^2_s(B_N)$ ($s>-1$), and obtained the following result concerning
compact difference of two special composition operators. About
compact  differences of   more general linear fractional composition
operators, they raised the question ($\ast$). \\ \\ {\bf Theorem
B.}\begin{em} Suppose $\varphi$ is a linear fractional self-map of
$B_N$ with a boundary fixed point. The operator
$C_{\varphi\circ\sigma}- C_{\sigma\circ\varphi}$ is compact on
$H^2(B_N)$ or $A^2_s(B_N)$ ($s>-1$) if and only if
$\varphi\circ\sigma=\sigma\circ\varphi$.\end{em} \\ \par In this
section, for linear fractional self-maps $\varphi$  and $\psi$ of
$B_N$, we will give a necessary and sufficient condition for
$C_\varphi-C_\psi$ to be compact on $H^2(B_N)$ or   $A^2_s(B_N)$
($s>-1$), which completely  answers  the question on compact
differences  of linear fractional composition operators in several
variables. In the proof of our  result,  Corollary 2.2 in Section 2
is   an essential tool. On the other hand, in order to  eventually
deduce that the symbols of two composition operators  are
equivalent, no matter  what  case in the disc or in the ball,
Bourdon \cite{Bo} and MacCluer and Weir \cite{MW} all used the fact
that if $C_\varphi-C_\psi$  is compact then
$$\rho(z_n)\frac{1-|z_n|^2}{1-|\varphi(z_n)|^2}$$ converges to zero for any sequence
$\{z_n\}$ with $|z_n|\to 1$.  This result will unavoidably be used
in our proof and so some of our treatments may be similar to those
in the proof of Theorem B. \\ \\ {\bf Theorem 3.1.}\begin{em}
Suppose that $\varphi$ and $\psi$ are linear fractional self-maps of
$B_N$. The difference $C_{\varphi}-C_{\psi}$ is compact on
$H^2(B_N)$ or $A^2_s(B_N)$ ($s>-1$) if and only if either both
$C_{\varphi}$ and $C_{\psi}$ are compact or $\varphi=\psi$.\end{em}
\\ \\  Proof. The sufficient condition is trivial, so we only need
to prove the necessity. If $C_{\varphi}-C_{\psi}$ is compact, it is
easy to see that $C_{\varphi}$ and $C_{\psi}$  must be compact or
not at the same time. Now, we assume that $C_\varphi$ is not
compact, then there exist $\zeta$ and $\eta$  on $\partial B_N$ such
that $\varphi(\zeta)=\eta$ (here, we have used the fact that
$C_{\varphi}$  is compact  if and only if $||\varphi||_{\infty}<1$
 for linear fractional self-map $\varphi$ of $B_N$). It follows that
$\varphi$ has finite angular derivative at $\zeta$ from the
smoothness of $\varphi$ on $\overline{B_N}$. Write
$t=d_\varphi(\zeta)$, thus $0<t<\infty$ by the
Julia-Carath\'{e}odory Theorem in $B_N$. Moreover, applying
Corollary 2.2, we  see that $\psi(\zeta)=\varphi(\zeta)$ and
$d_{\psi}(\zeta)=d_{\varphi}(\zeta)$  must hold. Similar to the
proof of Theorem 2.1, it suffices to
 assume that $\zeta=\eta= e_1$.

First, we give some information for the  matrix
\begin{displaymath} m_{\varphi}=\left(\begin{array}{ccc}
 A & B \\ C^{\ast} & d \end{array}\right)\end{displaymath} associated with $\varphi$,
where $A=(a_{jk})$,  $B=(b_j)$ and $C=(c_j)$ for $j, k=1, \ldots,
N$, and  $d>0$.  Since $\varphi(e_1)=e_1$, this  gives
$$a_{11}+b_1=\overline{c_1}+d\eqno(3.1)$$ and $a_{j1}+b_j=0$ for
$2\le j \le N$.  Note that Lemma 6.6 of \cite{CM1} shows  that $D_j
\varphi_1(e_1)=0$  and by
 Equation $(3.1)$, we compute that  $D_j
\varphi_1(e_1)=(a_{1j}-\overline{c_j})/(\overline{c_1}+d)$ for $j=2,
\ldots, N$. Thus $a_{1j}=\overline{c_j}$    for  $j=2,  \ldots, N$.
On the other hand, by the Julia-Carath\'{e}odory Theorem in $B_N$,
we   have $t=d_{\varphi}(e_1)=D_1
\varphi_1(e_1)=(a_{11}-\overline{c_1})/(\overline{c_1}+d)$.  The
arguments above yield  that
\begin{displaymath}
m_\varphi=\left(
\begin{array}{ccccc}
 a_{11} & \overline{c_2} & \ldots & \overline{c_N} & \overline{c_1}+d-a_{11}
\\ -b_2 &  &  &  & b_2 \\ \vdots & & & & \vdots \\ -b_N & & & & b_N
\\ \overline{c_1} & \overline{c_2} & \ldots & \overline{c_N} & d
\end{array}\right).
\end{displaymath}
Since  $\overline{c_1}+d\ne 0$,  setting
$K=\overline{c_1}/(\overline{c_1}+d)$,
$\beta_j=b_j/(\overline{c_1}+d)$  and
$\gamma_j=\overline{c_j}/(\overline{c_1}+d)$ for $2\le j\le N$, we
obtain an equivalent matrix for $\varphi$ up to multiply all entries
of $m_\varphi$ by the value $(\overline{c_1}+d)^{-1}$, that is
\begin{displaymath} \varphi\sim T\equiv
\left(\begin{array}{ccccc} t+K & \gamma_2 & \ldots & \gamma_N &
1-t-K
\\ -\beta_2 &  &  &  & \beta_2 \\ \vdots & & & & \vdots \\ -\beta_N & & & & \beta_N
\\ K & \gamma_2 & \ldots & \gamma_N &  1-K
\end{array}\right).
\end{displaymath}
Note that we have proved that $\psi(e_1)=\varphi(e_1)=e_1$ and
$d_\psi(e_1)=d_\varphi(e_1)=t$. Using  similar discussions as above,
we can get a matrix $S$ for $\psi$ with parameters $K'$, $\beta'_j$
and $\gamma'_j$ replacing $K$, $\beta_j$ and $\gamma_j$. For needed
later, we denote the $(j, k)$ entries of $T$ and $S$ by
$\alpha_{jk}$  and $\alpha'_{jk}$ respectively, where $j, k=2,
\ldots, N$.

First, we will prove $K=K'$ and  the argument is similar to Step 4
in the proof of Theorem B. For the convenience of the reader, we
will give a  proof in the case of the Hardy space. Define maps
$\rho$, $\rho': D\rightarrow D$
 by
$$\rho(\lambda)=\varphi_1(\lambda e_1)=\frac{(t+K)\lambda+1-t-K}{K\lambda+1-K},$$
 $$\rho'(\lambda)=\psi_1(\lambda e_1)=\frac{(t+K')\lambda+1-t-K'}{K'\lambda+1-K'}.$$
 If $K\ne K'$, it is clear that $\rho$ and $\rho'$ are distinct linear
fractional self-maps of $D$ with $\rho(1)=\rho'(1)=1$, then the
difference  $C_{\rho_j}-C_{\rho'_j}$  is not compact on
$A^2_{N-2}(D)$ (see \cite{Mo}). Hence, there exists a bounded
sequence $\{f_n\}$ in $A^2_{N-2}(D)$  which tends to $0$ uniformly
on compact subsets of $D$, and
$$||(C_{\rho_j}-C_{\rho'_j})f_n||_{A^2_{N-2}(D)}\nrightarrow 0$$ as
$n\to\infty$.  Define functions $F_n(z_1, z')=f_n(z_1)$  for $(z_1,
z')\in  B_N$, we  see that $\{F_n\}$  is a sequence tending to $0$
uniformly on compact subsets of $B_N$ and
$||F_n||_{H^2(B_N)}=||f_n||_{A^2_{N-2}(D)}$ (see 1.4.5 in
\cite{Ru}).  Setting   $g_n(\lambda)=F_n\circ \varphi(\lambda
e_1)-F_n\circ \psi(\lambda e_1)$ for $\lambda\in D$,  by Proposition
2.21 in \cite{CM1} this defines a restriction operator satisfying
$$||F_n\circ \varphi-F_n\circ \psi||_{H^2(B_N)}\ge ||g_n||_{A^2_{N-2}(D)}.$$
Obviously, we have $g_n(\lambda)=f_n\circ \varphi_1(\lambda
e_j)-f_n\circ \psi_1(\lambda
e_j)=f_n\circ\rho(\lambda)-f_n\circ\rho'(\lambda)$. Therefore,
$$||F_n\circ \varphi-F_n\circ \psi||_{H^2(B_N)}\ge
||f_n\circ\rho_j-f_n\circ\rho'_j||_{A^2_{N-2}(D)}$$  and so
$||(C_{\varphi}-C_{\psi})F_n||_{H^2(B_N)}$  does not converge to $0$
as $n\to\infty$. This contradicts with the fact that
$C_{\varphi}-C_{\psi}$ is compact  and so  $K=K'$ holds.

Now,  for $2\le k\le N$ and $2\le j\le N$, a computation shows that
$$D_1\varphi_1(e_1)=t, D_k\varphi_1(e_1)=0, $$
$$D_{11}\varphi_1(e_1)=-2tK,
D_{1k}\varphi_1(e_1)=-t\gamma_k,
 D_{kk}\varphi_1(e_1)=0,$$
and $$D_1\varphi_j(e_1)=-\beta_j, D_k\varphi_j(e_1)=\alpha_{jk},$$
 $$D_{11}\varphi_j(e_1)=2K\beta_j,  D_{1k}\varphi_j(e_1)=-K
\alpha_{jk}+\gamma_k\beta_j,
 D_{kk}\varphi_j(e_1)=-2\gamma_k\alpha_{jk},$$
where $\varphi_j$  denotes the $j^{\,th}$ component of $\varphi$.
Since $\varphi$ is holomorphic in a neighborhood of $e_1$, we have
the following expansions:
$$\varphi_1(z_1, 0')
 = 1+t(z_1-1)-tK(z_1-1)^2+o(|z_1-1|^2)\eqno(3.2)$$
and $$\varphi_j(z_1, 0')
 =-\beta_j(z_1-1)+K\beta_j(z_1-1)^2+o(|z_1-1|^2) \qquad   \mbox{as}\   (z_1, 0')\to e_1.\eqno(3.3)$$

Let $\Gamma\equiv \{z=(z_1, 0')\in B_N: 1-|z_1|^2=|1-z_1|^2\}$ be a
$e_1$-- curve.  It is clear that  $\mbox{Re}\, z_1=|z_1|^2$  for any
$z\in \Gamma$. Thus, for  points $z=(z_1, 0')\in \Gamma$, we use
(3.2) and (3.3) to obtain
{\setlength\arraycolsep{2pt}\begin{eqnarray*} 1-|\varphi(z)|^2 &=&
-2t\mbox{Re}\, (z_1-1)-t^2|z_1-1|^2+2t\mbox{Re}\,[K(z_1-1)^2]\\
&& -|z_1-1|^2\sum_{j=2}^N|\beta_j|^2+o(|z_1-1|^2)
\\ &=& (2t-t^2)|z_1-1|^2+2t\mbox{Re}\,[K(z_1-1)^2]-|z_1-1|^2\sum_{j=2}^N|\beta_j|^2+o(|z_1-1|^2)
.\hspace{9.4cm}
\end{eqnarray*}}Note that $(z_1-1)^2/|z_1-1|^2\to -1$  as $z$ approaches $e_1$ along
$\Gamma$, so $(1-|z|^2)/(1-|\varphi(z)|^2)$ has a finite limit
$$\frac{1}{2t-t^2-2t\mbox{Re}\, K- \sum_{j=2}^N |\beta_j|^2}.$$
However, combining the compactness of $C_\varphi-C_\psi$ with
Theorem 3 in \cite{MW}, we know that
$$\rho(z)\frac{1-|z|^2}{1-|\varphi(z)|^2}$$ must tend to $0$ as $z$ approaches the boundary of $B_N$.
This forces    that $\rho(z)$  goes to $0$ as $z$ approaches $e_1$
along $\Gamma$, so
$$1-\rho^2(z)=\frac{(1-|\varphi(z)|^2)(1-|\psi(z)|^2)}{|1-<\varphi(z),
\psi(z)>|^2}$$ converges  to $1$.

Using $K'=K$ and similar expansions as $(3.2)$ and $(3.3)$ for the
components of $\psi$ with $\beta_j$ replaced by $\beta'_j$, we
calculate that
$$\frac{1-|\psi(z)|^2}{1-|z|^2}\rightarrow 2t-t^2-2t\mbox{Re}\, K-
\sum_{j=2}^N |\beta'_j|^2$$ and {\setlength\arraycolsep{2pt}
\begin{eqnarray*} 1-<\varphi, \psi> &=& 1-|\psi|^2-<\varphi-\psi, \psi>
\\ &=&
1-|\psi|^2-|z_1-1|^2\sum_{j=2}^N\overline{\beta'_j}(\beta_j-\beta'_j)+o(|z_1-1|^2)
\end{eqnarray*}}as $z\to e_1$ along the curve $\Gamma$.
Write $a=t-t^2-2t\mbox{Re}\, K$,  because
$$t=\liminf\limits_{z\to
e_1}\frac{1-|\varphi(z)|^2}{1-|z|^2}=\liminf\limits_{z\to
e_1}\frac{1-|\psi(z)|^2}{1-|z|^2}$$ as $z$  approaches $e_1$
unrestrictedly,  this implies that $a-\sum_{j=2}^N |\beta_j|^2\ge 0$
and $a-\sum_{j=2}^N |\beta'_j|^2\ge 0$. Therefore,  as $z$ tends to
$e_1$  along $\Gamma$, $1-\rho^2(z)$ converges to
$$\frac{(t+a-\sum_{j=2}^N
|\beta_j|^2)(t+a- \sum_{j=2}^N |\beta'_j|^2)}{|t+a-\sum_{j=2}^N
\beta_j\overline{\beta'_j}|^2}.\eqno(3.4)$$

Set $L_1=(\sum_{j=2}^N |\beta_j|^2)^{1/2}$, $L_2=(\sum_{j=2}^N
|\beta'_j|^2)^{1/2}$ and $I=\sum_{j=2}^N
\beta_j\overline{\beta'_j}$. Since $|I|\le L_1L_2$, one may write
$|I|=\theta L_1L_2$ for $\theta\in [0, 1]$.   The  above discussion
gives that  the value in (3.4) is equal to $1$, it follows that  $$
0 = \biggl|t+a-\sum_{j=2}^N
\beta_j\overline{\beta'_j}\biggr|^2-\biggl(t+a-\sum_{j=2}^N
|\beta_j|^2\biggr)\biggl(t+a- \sum_{j=2}^N |\beta'_j|^2\biggr) $$ $$
=  (t+a)(L_1^2+L_2^2-2\mbox{Re}\, I)+|I|^2-L_1^2L_2^2 $$  $$  \ge
(t+a)(L_1^2+L_2^2-2|I|)+|I|^2-L_1^2L_2^2 \eqno(3.5) $$ $$ =
(t+a)(L_1^2+L_2^2-2\theta L_1L_2)+(\theta^2-1)L_1^2L_2^2.$$ Let
$h(\theta)=(t+a)(L_1^2+L_2^2-2\theta L_1L_2)+(\theta^2-1)L_1^2L_2^2$
for fixed $L_1$ and $L_2$. Observe that the derivative
$h'(\theta)=-2(t+a)L_1L_2+2\theta L_1^2L_2^2= 2L_1L_2(\theta
L_1L_2-a-t)<0$  from $a\ge L_1^2$  and $a\ge L_2^2$. Thus, the
function $h(\theta)$  is  decreasing on $[0, 1]$  and $0\ge
h(\theta)\ge h(1)$ for any $\theta\in [0, 1]$. However,
$h(1)=(t+a)(L_1-L_2)^2\ge 0$, this gives  $h(1)=0$ and then
$L_1=L_2$. Hence, the equalities in (3.5) and  $h(\theta)\ge h(1)$
must be attained,  which implies that $\mbox{Re}\, I=|I|$  and
$\theta=1$. All these  force  that   $\beta_j=\beta'_j$ for $2\le
j\le N$.

Next, we will prove that $\alpha_{jk}=\alpha'_{jk}$  for $j, k=2,
\ldots, N$. Fix $k\ge 2$, let $\Gamma_k(r)=re_1+\sqrt{1-r}e_k$ for
$0<r<1$. As $r\to 1^-$, the following expansions hold,
$$\varphi_1(r e_1+\sqrt{1-r}
e_k)=1+t(r-1)+o(1-r)$$   and  for $j\ge 2$, $$\varphi_j(r
e_1+\sqrt{1-r}
e_k)=-\beta_j(r-1)+\alpha_{jk}\sqrt{1-r}-\gamma_k\alpha_{jk}(1-r)+o(1-r).$$
Since   $\beta_j=\beta'_j$  for $2\le j \le N$, we can obtain
analogous expansions for the coordinates of $\psi$ with parameters
$\gamma_k, \alpha_{jk}$ replaced by
 $\gamma'_k, \alpha'_{jk}$. A computation shows that  $$\frac{1-|z|^2}{1-|\varphi(z)|^2}\rightarrow
\frac{1}{2t-\sum^N_{j=2} |\alpha_{jk}|^2}$$ and
$$1-\rho^2(z)\rightarrow\frac{(2t-\sum^N_{j=2}
|\alpha_{jk}|^2)(2t-\sum^N_{j=2} |\alpha'_{jk}|^2)}{|2t-\sum^N_{j=2}
\alpha_{jk}\overline{\alpha'_{jk}}|^2}$$ as $z$ tends to $e_1$ along
the curve $\Gamma_k$. Similarly,  using
$t=d_\varphi(e_1)=d_\psi(e_1)$, we  deduce that $\sum^N_{j=2}
|\alpha_{jk}|^2\le t$ and $\sum^N_{j=2} |\alpha'_{jk}|^2\le t$. Let
$L_1=(\sum^N_{j=2} |\alpha_{jk}|^2)^{1/2}$,  $L_2=(\sum^N_{j=2}
|\alpha'_{jk}|^2)^{1/2}$  and $I=\sum^N_{j=2}
\alpha_{jk}\overline{\alpha'_{jk}}$. Since $C_\varphi-C_\psi$ is
compact,  by Theorem 3 of \cite{MW}, we see that $1-\rho^2(z)$ must
tend to $1$ as $z\to e_1$ along $\Gamma_k$ and so
{\setlength\arraycolsep{2pt}
\begin{eqnarray*}0&=&\biggl|2t-\sum^N_{j=2}
\alpha_{jk}\overline{\alpha'_{jk}}\biggr|^2-\biggl(2t-\sum^N_{j=2}
|\alpha_{jk}|^2\biggr)\biggl(2t-\sum^N_{j=2} |\alpha'_{jk}|^2\biggr)
\\ &=& 2t(L_1^2+L_2^2-2\mbox{Re}\, I)+|I|^2-L_1^2L_2^2.
\end{eqnarray*}}The  remaining proof is similar to  that used in the proof of   $\beta_j=\beta'_j$,   so we omit it.

Finally, it remains to prove $\gamma_k=\gamma'_k$ for $ 2\le k\le
N$. Fix $k\ge 2$  and $0<r<1/2$, let $\Gamma_{k, r}
\equiv\{z=z_1e_1+(z_1-1)e_k\in B_N, z_1=1-r+re^{i\theta}\ \mbox{for
\ real }\  \theta\}$. Then  the curve $\Gamma_{k, r}$ approaches
$e_1$  as $\theta\to 0$, i.e as $z\to e_1$ and  for points $z\in
\Gamma_{k, r}$, we have $\frac{1-|z_1|^2}{|1-z_1|^2}=\frac{1-r}{r}$,
$\frac{1-|z|^2}{|1-z_1|^2}=\frac{1-2r}{r}$ and
$\frac{1-Re\,z_1}{|1-z_1|^2}=\frac{1}{2r}$. As $z$ tends to $e_1$
along $\Gamma_{k, r}$, we get that
$$\varphi_1(z)
 = 1+t(z_1-1)-t(K+\gamma_k)(z_1-1)^2+o(|z_1-1|^2),$$
$$\varphi_j(z)
 =-(\beta_j-\alpha_{jk})(z_1-1)+(K+\gamma_k)(\beta_j-\alpha_{jk})(z_1-1)^2+o(|z_1-1|^2)$$
for $j\ge 2$ and  similar expansions for the components of $\psi$
only with $\gamma'_k$ replacing $\gamma_k$. Therefore,
$$\frac{1-|z|^2}{1-|\varphi(z)|^2}\rightarrow
\frac{1-2r}{t-r[t^2+2t\mbox{Re}\,(K+\gamma_k)+\sum^N_{j=2}
|\beta_j-\alpha_{jk}|^2]}$$  as $z$ approaches   $e_1$ along
$\Gamma_{k, r}$. By Theorem 3 in \cite{MW}, this with the
compactness of $C_{\varphi}- C_{\psi}$
 shows that $\rho(z)$ converges to
zero as $z\to e_1$ along $\Gamma_{k, r}$.  It is clear that
$$\rho(z)\ge \biggl|\frac{\varphi(z)-\frac{<\psi(z), \varphi(z)>}
{|\varphi(z)|^2}\varphi(z)}{1-<\psi(z), \varphi(z)>}\biggr|$$ and
{\setlength\arraycolsep{2pt}
\begin{eqnarray*}\biggl|\frac{\varphi-\frac{<\psi, \varphi>}{|\varphi|^2}\varphi}{1-<\psi, \varphi>}\biggr|
&=&\frac{1}{|\varphi|}\biggl|\frac{<\varphi-\psi,
\varphi>}{1-|\varphi|^2+<\varphi-\psi, \varphi>}\biggr|
\\ &   =& \frac{1}{|\varphi|}\biggl|\frac{-t(\gamma_k-\gamma'_k)(z_1-1)^2+o(|z_1-1|^2)}
{1-|\varphi|^2-t(\gamma_k-\gamma'_k)(z_1-1)^2+o(|z_1-1|^2)}\biggr|
\\ & =& \frac{1}{|\varphi|}\biggl|\frac{-t(\gamma_k-\gamma'_k)\frac{(z_1-1)^2}{|z_1-1|^2}+o(1)}
{\frac{1-|\varphi|^2}{|z_1-1|^2}-t(\gamma_k-\gamma'_k)\frac{(z_1-1)^2}{|z_1-1|^2}+o(1)}\biggr|.
\end{eqnarray*}}Taking  the limit  as $z\to e_1$  along $\Gamma_{k, r}$
and using $\frac{(z_1-1)^2}{|z_1-1|^2} \to -1$, the quotient
converges to
$$\biggl|\frac{t(\gamma_k-\gamma'_k)}{\frac{t}{r}-t^2-2t\mbox{Re}\,K-t(\overline{\gamma_k}+\gamma'_k)-\sum^N_{j=2}
|\beta_j-\alpha_{jk}|^2}\biggr|.$$ It must be  zero from the above
discussions.  So we have  $\gamma_k=\gamma'_k$ as desired   and
complete the proof. \ \ $\Box$

Now, Combining the argument in Section 1 with the proof of Theorem B
or as a corollary of   Theorem 3.1, Theorem B can be improved as
follows. \\ \\ {\bf Theorem 3.2.}\begin{em} Suppose  that $\varphi$
is a linear fractional self-map of $B_N$ with
$||\varphi||_\infty=1$. The operator $C_{\varphi\circ\sigma}-
C_{\sigma\circ\varphi}$ is compact on $H^2(B_N)$ or $A^2_s(B_N)$
($s>-1$) if and only if
$\varphi\circ\sigma=\sigma\circ\varphi$.\end{em} \\ \\   Proof. We
only need to prove one direction. If $||\varphi||_\infty=1$, there
exist $\zeta$ and $\eta$ on $\partial B_N$ such that
$\varphi(\zeta)=\eta$, then $\sigma(\eta)=\zeta$ by Lemma 1 of
\cite{MW}. This implies that $||\sigma\circ\varphi||_\infty=1$ for
the linear fractional map $\sigma\circ\varphi$ of $B_N$ and so
$C_{\sigma\circ\varphi}$ is not compact. Therefore, if
$C_{\varphi\circ\sigma}- C_{\sigma\circ\varphi}$ is compact, by
Theorem 3.1,  we have $\varphi\circ\sigma=\sigma\circ\varphi$.  \ \
$\Box$

\small DEPARTMENT OF MATHEMATICS, TONGJI UNIVERSITY,

SHANGHAI 200092, CHINA

DEPARTMENT OF APPLIED  MATHEMATICS,

SHANGHAI FINANCE UNIVERSITY,

SHANGHAI  201209, CHINA

{\it E-mail address:} liangying1231@163.com \\ \par WUHAN INSTITUTE
OF PHYSICS AND MATHEMATICS,

CHINESE ACADEMY OF SCIENCES, WUHAN 430071, CHINA

{\it E-mail address:}  ouyang@wipm.ac.cn

\end{document}